\documentclass[letterpaper, 10 pt, conference]{ieeeconf}  

\IEEEoverridecommandlockouts                              
\usepackage[utf8]{inputenc}
\usepackage[T1]{fontenc}
\usepackage[english]{babel} 

\usepackage{ae}
\usepackage{aecompl}
\usepackage{aeguill}

\usepackage{graphicx} 

\usepackage{amsfonts,amsmath,amssymb} 
\usepackage{array}
\usepackage{subfigure}
\usepackage{cite}

\graphicspath{{Figures/}}
\usepackage{tikz}
\usepackage{pgfplots}
\pgfplotsset{compat=newest}
\pgfplotsset{plot coordinates/math parser=false}
\tikzstyle{dotted} = [line width = 0.5pt, dash pattern = on \pgflinewidth off 2\pgflinewidth]
\pgfplotsset{%
	tick label style={font=\scriptsize,%
		/pgf/number format/precision=3,%
		/pgf/number format/std=-3:3},%
	legend style={font=\small},%
	label style={font=\small},%
	every axis post/.append style={label style={font=\small},%
		axis line style=-,%
		scaled ticks=false}%
}

\newlength{\figurewidth}
\newlength{\figureheight}

\providecommand{\norm}[1]{\left\lVert#1\right\rVert}

\def\ordfield{}
\def\subfield{}
\def\expfield{}
\def\getfields#1{%
	\def\ordfield{}%
	\def\subfield{}%
	\def\expfield{}%
	\expandafter\getexpfield\expandafter\beginsym#1^!\endsymb%
}
\def\extractsubfield\beginsymb#1_!\endsymb{%
	\def\subfield{#1}%
}
\def\getsubfield\beginsymb#1_#2\endsymb{%
	\if!#2\relax\else%
		\extractsubfield\beginsymb#2\endsymb%
	\fi%
	\def\ordfield{#1}%
}
\def\extractexpfield\beginsymb#1^!\endsymb{%
	\getsubfield\beginsymb#1_!\endsymb%
	\expandafter\def\expandafter\expfield\expandafter{\ordfield}
}
\def\getexpfield\beginsym#1^#2\endsymb{%
	\if!#2\relax\else%
		\extractexpfield\beginsymb#2\endsymb%
	\fi%
	\getsubfield\beginsymb#1_!\endsymb%
}

\newcommand{\be}{\begin{equation}}
\newcommand{\ee}{\end{equation}}
\newcommand{\beno}{\begin{equation*}}
\newcommand{\eeno}{\end{equation*}}
\newcommand{\barCl}{\begin{IEEEeqnarray}{rCl}}
\newcommand{\earCl}{\end{IEEEeqnarray}}
\newcommand{\barClno}{\begin{IEEEeqnarray*}{rCl}}
\newcommand{\earClno}{\end{IEEEeqnarray*}}
\newcommand{\ans}{\IEEEeqnarraynumspace}
\newcommand{\ase}{\IEEEyessubnumber}
\newcommand{\ane}{\IEEEnonumber}

\newcommand{\AB}{{\alpha\beta}}		
\newcommand{\dq}{{dq}}				
\newcommand{\Hdq}{\mathcal{H}^\dq}	
\newcommand{\Hmdq}{\Hdq_m}			
\newcommand{\hmdq}{H_m^\dq}			
\newcommand{\usdq}{u_s^\dq}			
\newcommand{\usAB}{u_s^\AB}	
\newcommand{\isAB}{\imath_s^\AB}	
\newcommand{\isdq}{\imath_s^\dq}	
\newcommand{\irdq}{\imath_r^\dq}	
\newcommand{\phisdq}{\phi_s^\dq}	
\newcommand{\phirdq}{\phi_r^\dq}	
\newcommand{\Te}{T_e}				
\newcommand{\Tl}{T_l}				
\newcommand{\km}{\rho}				
\newcommand{\w}{\omega}				
\newcommand{\ws}{\w_s}				
\newcommand{\wg}{\w_g}				
\newcommand{\ts}{\theta_s}			
\newcommand{\ti}{\theta_i}			

\newcommand{\Lm}{L_m}				
\newcommand{\Ll}{L_l}				
\newcommand{\epsm}{\varepsilon_m}	
\newcommand{\epsl}{\varepsilon_l}	
\newcommand{\Rs}{R_s}				
\newcommand{\Rr}{R_r}				
\newcommand{\Jl}{J_l}				
\newcommand{\np}{n}					
\newcommand{\Wr}{\Omega_r}
\newcommand{\Xr}{X_r}
\newcommand{\W}{\Omega}				
\newcommand{\finj}{s}
\newcommand{\Finj}{S}
\newcommand{\Sal}{\mathbf{S}}		
\newcommand{\Id}{\mathcal{I}}		
\newcommand{\J}{\mathcal{J}}		
\newcommand{\Rot}{\mathcal{R}}		
\newcommand{\Sym}{\mathcal{S}}		
\newcommand{\Z}[2]{0_{#1,#2}}		

\newcommand{\transpose}[1]{{#1}^T}
\newcommand{\Pderive}[2]{\frac{\partial#1}{\partial#2}}
\newcommand{\pderive}[2]{\partial_{#2}#1}
\newcommand{\Tderive}[2][1]{%
	\if1#1\relax
		\frac{d#2}{dt}%
	\else%
		\frac{d^{#1}#2}{dt^{#1}}%
	\fi%
}
\newcommand{\tderive}[2][1]{
	\if1#1\relax
		\dot{#2}%
	\else%
		\if2#1\relax%
			\ddot{#2}%
		\else%
			{#2}^{(#1)}%
		\fi%
	\fi%
}
\newcommand{\lc}{LC}

\DeclareMathOperator{\Landau}{\mathcal{O}}

\newcommand{\lf}[1]{\overline{#1}}
\newcommand{\hf}[1]{\widetilde{#1}}

\newcommand{\equ}[1]{%
	\getfields{#1}%
	\ifx\subfield\empty%
		\ordfield_{eq}^{\expfield}%
	\else%
		\ordfield_{\subfield,eq}^{\expfield}%
	\fi%
}
\newcommand{\Equ}[1]{\equ{#1}}

\newcommand{\nbdash}{\nobreakdash-\hspace{0pt}}

\title{\LARGE \bf An analysis of the benefits of signal injection for low-speed sensorless control of induction motors}
\author{Pascal Combes\textsuperscript{1}, François Malrait\textsuperscript{2}, Philippe Martin\textsuperscript{1} and Pierre Rouchon\textsuperscript{1}
\thanks{\textsuperscript{1}~P.~Combes, P.~Martin and P.~Rouchon are with the Centre Automatique et Systèmes, MINES ParisTech, PSL Research University, Paris, France
{\tt\footnotesize\{pascal.combes, philippe.martin,pierre.rouchon\}@mines-paristech.fr}}%
\thanks{\textsuperscript{2}~F.~Malrait is with Schneider Toshiba Inverter Europe, 27120~Pacy-sur-Eure,~France
{\tt\footnotesize francois.malrait@schneider-electric.com}}
}

\begin{document}

\maketitle
\thispagestyle{empty}
\pagestyle{empty}

\begin{abstract}
	We analyze why low-speed sensorless control of the IM is intrinsically difficult, and what is gained by signal injection. The explanation relies on the control-theoretic concept of observability applied to a general model of the saturated IM. We show that the IM is not observable when the stator speed is zero in the absence of signal injection, but that observability is restored thanks to signal injection and magnetic saturation. The analysis also reveals that existing sensorless algorithms based on signal injection may perform poorly for some IMs under particular operating conditions. The approach is
	illustrated by simulations and experimental data.
\end{abstract}



\section{Introduction}
	``Sensorless'' control algorithms have become a standard feature of Variable Speed Drives for AC motors (``sensorless'' meaning that the algorithm is able to control the motor with only stator current measurements). While sensorless control at medium to high speed is now rather well understood, it remains a challenge at low speed, in particular for the Induction Motor (IM), see e.g.~\cite{LascuBB2005ITOIA,BenlaDCO2015ITOEC}.

A promising technique to tackle the problem is High-Frequency (HF) signal injection, originally proposed by~\cite{JanseL1995IToIA} in 1995. Since then, quite a lot of algorithms have been proposed for the Permanent Magnet Synchronous Motor, but much fewer for the IM.
Part of the problem is that most IMs do not exhibit any geometric saliency, so signal-injection-based algorithms have to rely on magnetic saliency~\cite{YoonS2014ITOPE,Brandstetter2011,Zatocil2008}, except for rather specific IMs with open slots~\cite{HoltzP2004IToIA} or suitable modifications~\cite{DegneL2000IToIA}. What is lacking is a good control model of the saturated IM, together with a thorough analysis of the effects of signal injection. A convincing physical explanation is proposed in~\cite{HaSIMS2000CRotIIAC,YoonS2014ITOPE}, but it is not completely satisfying as it relies on several physical approximations.

In this paper, we try to explain why low-speed sensorless control of the IM is intrinsically difficult and what is gained by signal injection. This explanation is based on the concept of observability, which expresses the possibility of theoretically recovering all the state variables of the system from the known signals and their derivatives (in our case, recovering the stator flux, rotor flux, rotor velocity and load torque from the measured stator currents and the input stator voltages). The concrete use of this property is the following: around a point where the system is not observable, it is very difficult to design a control law, even if the system were perfectly known; on the contrary, around a point where the system is (first-order) observable, it is always possible to design at least a basic control law relying for instance on an extended Kalman filter. We model the saturated IM following the approach of \cite{Jebai2014} and analyze the effects of signal injection with the notion of virtual measurements introduced in \cite{Combes2016}. We show that the IM, saturated or not, is not observable when the stator speed is zero (which corresponds to a straight line in the velocity-torque plane) in the absence of signal injection, but that observability is restored thanks to signal injection. The analysis also reveals that the existing algorithms~\cite{Brandstetter2011,JanseL1995IToIA,YoonS2014ITOPE,Zatocil2008} may perform poorly for some IMs. Indeed, these algorithms first estimate the rotor and stator fluxes by algebraically inverting the four relations provided by the two current measurements and the two virtual measurements obtained by signal injection; but in some operating conditions, the two virtual measurements convey nearly the same information, which renders this inversion impossible in practice.
 
The paper is organized as follows: in section~\ref{sec:nohf} we propose a general model for the saturated IM and study its observability; in section~\ref{sec:hf} we analyze the effects of signal injection on this model and illustrate the theory by simulations and experimental data; finally in section~\ref{sec:obshf} we show that observability is recovered everywhere thanks to signal injection.

\section{The saturated IM without signal injection}\label{sec:nohf}
	\subsection{A general model of the saturated IM}\label{sec:model}
We use the energy-based modeling approach developed in~\cite{Jebai2014}, which is an application of classical analytical mechanics. This approach encodes the information in a single scalar function, and automatically provides nonlinear models with a sensible physical structure:
\begin{itemize}
	\item it gives an expression for the electromagnetic torque valid in the presence of saturation
	\item it justifies the modeling of saturation in the fictitious $dq$~frame. This is usually taken for granted, though saturation physically takes place in the $abc$~frame
	\item it automatically enforces the reciprocity conditions (see e.g.~\cite{MelkeW1990IAITo,Sauer1992ECITo}).
\end{itemize}
Moreover this approach does not rely on the motor internal layout, but only on basic geometric symmetries enjoyed by any well-built electric motor.


Applied to the IM, whose state variables are the vector of stator flux linkages~$\phisdq$, the vector of rotor flux linkages~$\phirdq$, the kinetic momentum~$\rho$ and the (electrical) rotor angle~$\theta$, this approach yields
\barCl\label{eqn:model:dynamic}
	\Tderive{\phisdq} &=& \usdq - \Rs\isdq - \J\ws\phisdq \ase \label{eqn:model:dynamic:stator} \\
	\Tderive{\phirdq} &=& -\Rr\irdq - \J(\ws - \w)\phirdq \ase \label{eqn:model:dynamic:rotor} \\
	\np\Tderive{\km} &=& \Te - \Tl \ase \label{eqn:model:dynamic:speed} \\
	\Tderive{\theta} &=& \w \ase, \label{eqn:model:dynamic:angle}
\earCl
where $\ws = \tderive{\ts}$~is the speed of rotation of the synchronous $dq$~frame and $\J$~is the matrix
\beno
\J := \begin{pmatrix}
	0 & -1 \\
	1 & 0
\end{pmatrix};
\eeno
$\Rs$~and $\Rr$~are the stator and rotor resistances and $\np$~is the number of pole pairs.
The $\dq$~frame is obtained from the stator $\AB$~frame by a known rotation of angle $\ts$; consequently $\dq$~and $\AB$~stator variables are related by $x_s^\dq := \Rot(-\ts)x_s^\AB$ whereas the relation between $\dq$~and $\AB$~rotor variables is $x_r^\dq := \Rot(-\ts + \theta)x_r^\AB$, where
\beno
	\Rot(\eta) := \begin{pmatrix}
		\cos{\eta} & -\sin{\eta} \\
		\sin{\eta} & \cos{\eta}
	\end{pmatrix}
\eeno
is the rotation of angle~$\eta$.
The control inputs are the vector of stator voltages~$\usdq$ ($\usdq$~and $\usAB$~are equivalent since $\ts$ is known) and~$\ws$. The load torque~$\Tl$ is an unmeasured disturbance.


The model~\eqref{eqn:model:dynamic} is then closed thanks to the energy function~$\Hdq$ (yet to be specified) according to the structural relations
\barCl\label{eqn:model:relations}
	\isdq &:=& \Pderive{\Hdq}{\phisdq}(\theta, \rho, \phisdq, \phirdq) \ase \label{eqn:model:isdq} \\
	\irdq &:=& \Pderive{\Hdq}{\phirdq}(\theta, \rho, \phisdq, \phirdq) \ase \label{eqn:model:irdq} \\
	\Te &:=& -\np\Pderive{\Hdq}{\theta} + \np\transpose{\phirdq} \J \irdq \ase \label{eqn:model:torque} \\
	\w &:=& \Pderive{\Hdq}{\km}(\theta, \km, \phisdq, \phirdq); \ase \label{eqn:model:speed}
\earCl
these equations relate the stator currents~$\isdq$, the rotor currents~$\irdq$, the (electrical) speed~$\w$ and the electromagnetic torque~$\Te$ to the state variables of~\eqref{eqn:model:dynamic}.

We now particularize the form of the energy function~$\Hdq$. On the one hand, $\Hdq$~is the sum of the  magnetic energy $\Hmdq$ and mechanical energy; as we are not interested here in mechanical modeling, we consider the simple mechanical energy~$\frac{\np^2}{2\Jl}\km^2$, which corresponds to a balanced inertial load with moment of inertia~$\Jl$.
As for the magnetic energy~$\Hmdq$, it obviously does not depend on the mechanical variable~$\km$;
moreover the dependence on the rotor angle~$\theta$ is small enough to be neglected on a well-built IM.
As a consequence the energy function~$\Hdq$ reads
\beno
	\Hdq(\theta, \km, \phisdq, \phirdq) = \frac{\np^2}{2\Jl}\km^2 + \Hmdq(\phisdq, \phirdq),
\eeno
and the structural relations~\eqref{eqn:model:relations} simplify to
\barCl
	\isdq &=& \Pderive{\Hmdq}{\phisdq}(\phisdq, \phirdq) \ase \label{eqn:sinus:is} \\
	\irdq &=& \Pderive{\Hmdq}{\phirdq}(\phisdq, \phirdq) \ase \label{eqn:sinus:ir} \\
	\Te &=& \np \transpose{\phirdq} \J \irdq \ase \label{eqn:sinus:Te} \\
	\w &=& \frac{\np^2}{\Jl}\km. \ase \label{eqn:sinus:omega}
\earCl
Moreover the state variable~$\theta$ can be dropped from the model \eqref{eqn:model:dynamic} since its defining equation~\eqref{eqn:model:dynamic:angle} is completely autonomous (cyclic variable).

On the other hand, we can further simplify~$\Hmdq$ by considering the construction symmetries of the electric motor. Indeed, IMs are built so that no privileged flux paths exist in the rotor and the stator, i.e. they are geometrically non-salient. This implies that rotating the stator and rotor fluxes by the same angle~$\eta$ preserves the energy, in other words,
\be\label{eqn:inv:rot}
	\Hmdq(\phisdq, \phirdq) = \Hmdq(\Rot(\eta)\phisdq, \Rot(\eta)\phirdq).
\ee
Besides, exchanging two phases on the stator and the rotor and symmetrizing the position of the rotor does not change the energy either, thus
\be\label{eqn:inv:sym}
	\Hmdq(\phisdq, \phirdq) = \Hmdq(\Sym\phisdq, \Sym\phirdq),
\ee
where $\Sym$~denotes the reflexion matrix
\beno
	\Sym := \begin{pmatrix}
		1 & 0 \\
		0 & -1
	\end{pmatrix}.
\eeno
The only scalar quantities depending on~$\phisdq$ and~$\phirdq$ which respect the constraints~\eqref{eqn:inv:rot} and~\eqref{eqn:inv:sym} are the norms $\norm{\phisdq}^2$ and $\norm{\phirdq}^2$ and the scalar product~$\transpose{\phisdq}\!\phirdq$. Consequently, the magnetic energy~$\Hmdq$ can be written as a function~$\hmdq$ depending only on three scalar variables
\be\label{eqn:model:sinus}
	\Hmdq(\phisdq, \phirdq) = \hmdq\Biggl(\frac{\norm{\phisdq}^2}{2}, \transpose{\phisdq}\!\phirdq, \frac{\norm{\phirdq}^2}{2}\Biggr).
\ee	
The currents~\eqref{eqn:sinus:is}--\eqref{eqn:sinus:ir} then read
\barCl\label{eqn:inv:currents}
	\isdq &=& \pderive{\hmdq}{1}(\ldots)\phisdq + \pderive{\hmdq}{2}(\ldots)\phirdq \ans \ase \\
	\irdq &=& \pderive{\hmdq}{2}(\ldots)\phisdq + \pderive{\hmdq}{3}(\ldots)\phirdq, \ans \ase 
\earCl
where $\pderive{\hmdq}{n}$~is the first partial derivative of~$\hmdq$ with respect to the $n^{\rm th}$ variable and $(\ldots)$~stands for $\Bigl(\norm{\phisdq}^2\!\!\big/2, \transpose{\phisdq}\!\phirdq, \norm{\phirdq}^2\!\!\big/2\Bigr)$.
As a consequence, the electromagnetic torque~$\Te$ can be calculated with any of the formulae
\beno
	\Te = \np \transpose{\phirdq}\J\irdq = -\np \transpose{\phisdq}\J\isdq.
\eeno
The second alternative greatly simplifies the observability study hereafter.

To conclude, the general model of the saturated IM is given by 
\barCl\label{eqn:dynamic}
	\Tderive{\phisdq} &=& \usdq - \Rs\isdq - \J\ws\phisdq \ase \label{eqn:dynamic:stator} \\
	\Tderive{\phirdq} &=& -\Rr\irdq - \J(\ws - \w)\phirdq \ase \label{eqn:dynamic:rotor} \\
	\frac{\Jl}{\np}\Tderive{\w} &=& -\np \transpose{\phisdq}\J\isdq - \Tl, \ase \label{eqn:dynamic:speed}
\earCl
with~$\isdq$ and~$\irdq$ given by~\eqref{eqn:inv:currents}.

The simplest choice for the energy function is the quadratic form
\be\label{eqn:model:linear}
	\Hmdq(\phisdq, \phirdq) = \frac{\norm{\phisdq + \phirdq}^2}{4(2\Lm + \Ll)} + \frac{\norm{\phisdq - \phirdq}^2}{4\Ll}, 
\ee
which corresponds to an unsaturated motor with linear current-flux relations. Notice~\eqref{eqn:model:linear} has indeed the form~\eqref{eqn:model:sinus}.
	\subsection{Sensorless first-order observability of the saturated IM}\label{sec:obsnohf}
Sensorless control of IMs is difficult at low speed because the equations are unobservable at zero stator speed. In other words, it is theoretically impossible to retrieve the state variables of~\eqref{eqn:dynamic} from the measurements~$\isAB$ and their time derivatives (or equivalently~$\isdq$, since the angle~$\ts$ between the~$\AB$ and~$\dq$ frames is known), as is well-known for the unsaturated model, see e.g~\cite{Glumineau2015}. In this section we generalize this result to the saturated IM model~\eqref{eqn:dynamic}. 
For the sake of simplicity, we restrict to the analysis of first-order observability, but the result can be generalized to nonlinear observability.


An equilibrium point of the dynamical system~\eqref{eqn:model:sinus} $(\Equ{\phisdq}, \Equ{\phirdq}, \equ{\w};\allowbreak \Equ{\usdq}, \equ{\ws}, \equ{\Tl})$ is defined by
\barClno
	0 &=& \Equ{\usdq} - \Rs\Equ{\isdq} - \J\equ{\ws}\Equ{\phisdq} \\
	0 &=& -\Rr\Equ{\irdq} - \J\equ{\wg}\Equ{\phirdq} \\
	0 &=& -\np \transpose{\Equ{\phisdq}}\J\Equ{\isdq} - \equ{\Tl},
\earClno
where $\equ{\wg} := \equ{\ws} - \equ{\w}$. The first-order linearization of~\eqref{eqn:dynamic} around this equilibrium reads
\barCl\label{eqn:model:linearized}
	\Tderive{\delta\phisdq} &=& \delta\usdq - \Rs\delta\isdq - \J\equ{\ws}\delta\phisdq - \J\Equ{\phisdq}\delta\ws \ans \ase \\
	\Tderive{\delta\phirdq} &=& -\Rr\delta\irdq - \J\equ{\wg}\delta\phirdq \ane \\
	& & -\>\J\Equ{\phirdq}(\delta\ws - \delta\w) \ans \ase \\
	\frac{\Jl}{\np}\Tderive{\delta\w} &=& \np\transpose{\Equ{\isdq}}\J\delta\phisdq - \np\transpose{\Equ{\phisdq}}\J\delta\isdq - \delta\Tl, \ans \ase 
\earCl
with the linearized current-flux relations
\barCl\label{eqn:obsnohf:currents}
	\delta\isdq = \frac{\partial^2\Hmdq}{\partial{\phisdq}^2}\delta\phisdq + \frac{\partial^2\Hmdq}{\partial{\phisdq}\partial{\phirdq}}\delta\phirdq \ase \label{eqn:obsnohf:is} \\
	\delta\irdq = \frac{\partial^2\Hmdq}{\partial{\phirdq}\partial{\phisdq}}\delta\phisdq + \frac{\partial^2\Hmdq}{\partial{\phirdq}^2}\delta\phirdq. \ase \label{eqn:obsnohf:ir}
\earCl
In~\eqref{eqn:obsnohf:currents} and in the rest of this section, the partial derivatives of~$\Hmdq$ are evaluated at~$(\Equ{\phisdq}, \Equ{\phirdq})$, which is not explicitly written for lack of space. Notice the 2x2 matrices $\frac{\partial^2\Hmdq}{\partial{\phisdq}^2}$, $\frac{\partial^2\Hmdq}{\partial{\phisdq}\partial{\phirdq}}$ and $\frac{\partial^2\Hmdq}{\partial{\phirdq}^2}$ are invertible, as the energy function must be non-degenerate to be physically acceptable.

To study the first-order observability, we compute the time derivatives of the measurements. The first one is
\barClno
	\Tderive{\delta\isdq} &=& \frac{\partial^2\Hmdq}{\partial{\phisdq}^2}\Tderive{\delta\phisdq} + \frac{\partial^2\Hmdq}{\partial{\phisdq}\partial{\phirdq}}\Tderive{\delta\phirdq} \\
	&=& \Bigl(-\Rr\frac{\partial^2\Hmdq}{\partial{\phisdq}\partial{\phirdq}}\frac{\partial^2\Hmdq}{\partial{\phirdq}\partial{\phisdq}}  - \frac{\partial^2\Hmdq}{\partial{\phisdq}^2}\J\ws\Bigr)\delta\phisdq \\
	& &+\>\frac{\partial^2\Hmdq}{\partial{\phisdq}\partial{\phirdq}}\Bigl(-\Rr\frac{\partial^2\Hmdq}{\partial{\phirdq}^2} - \J\equ{\wg}\Bigr)\delta\phirdq \\
	& &+\>\frac{\partial^2\Hmdq}{\partial{\phisdq}\partial{\phirdq}}\J\Equ{\phirdq}\delta\w + \lc(\delta\usdq, \delta\ws, \delta\isdq),
\earClno
where $\lc$~denotes a linear combination of its arguments. Using~\eqref{eqn:obsnohf:is}, $\delta\phirdq$~can be expressed as a function of~$\delta\phisdq$ and~$\delta\isdq$,
\beno
	\delta\phirdq = \Biggl(\frac{\partial^2\Hmdq}{\partial{\phisdq}\partial{\phirdq}}\Biggr)^{-1}\Biggl(\delta\isdq - \frac{\partial^2\Hmdq}{\partial{\phirdq}^2}\delta\phisdq\Biggr),
\eeno
and can be replaced in the expression of~$\Tderive{\delta\isdq}$.
After some algebra this yields
\barCl\label{eqn:obsnohf:1}
	\Tderive{\delta\isdq} &=& -M\delta\phisdq + \frac{\partial^2\Hmdq}{\partial{\phisdq}\partial{\phirdq}}\J\Equ{\phirdq}\delta\w \ane \\
			& & +\>\lc(\delta\usdq, \delta\ws, \delta\isdq),
\earCl
where
\barClno
	M &:=& \frac{\partial^2\Hmdq}{\partial{\phisdq}\partial{\phirdq}}\Wr + \frac{\partial^2\Hmdq}{\partial{\phisdq}^2}\J\equ{\ws} - \Xr\frac{\partial^2\Hmdq}{\partial{\phisdq}^2} \equ{\wg} \\
	\Wr &:=& \Rr \frac{\partial^2\Hmdq}{\partial{\phirdq}\partial{\phisdq}} - \Rr \frac{\partial^2\Hmdq}{\partial{\phirdq}^2}\Biggl(\frac{\partial^2\Hmdq}{\partial{\phisdq}\partial{\phirdq}}\Biggr)^{-1}\frac{\partial^2\Hmdq}{\partial{\phisdq}^2} \\
	\Xr &:=& \frac{\partial^2\Hmdq}{\partial{\phisdq}\partial{\phirdq}} \J \Biggl(\frac{\partial^2\Hmdq}{\partial{\phisdq}\partial{\phirdq}}\Biggr)^{-1}
\earClno
($M$, $\Wr$ and~$\Xr$ of course depends on $\Equ{\phisdq}$ and~$\Equ{\phirdq}$).
The matrix~$M$ is invertible for any reasonable operating point; more precisely, it is never very far from $\lambda\Id + \mu\J$, since the saturated model is never very far from the unsaturated one.

The second time derivative of~$\delta\isdq$ is
\barClno
	\Tderive[2]{\delta\isdq} &=& M\J\equ{\ws}\delta\phisdq + \frac{\np^2}{\Jl} \frac{\partial^2\Hmdq}{\partial{\phisdq}\partial{\phirdq}}\J\Equ{\phirdq}\transpose{\Equ{\isdq}}\!\!\!\J\delta\phisdq \ane \\
	& & -\>\frac{\np}{\Jl}\frac{\partial^2\Hmdq}{\partial{\phisdq}\partial{\phirdq}}\J\Equ{\phirdq}\delta\Tl \\
	& & +\>\lc(\delta\usdq, \delta\ws, \delta\isdq, \tderive{\delta\isdq}).
\earClno
Adding $M \J M^{-1}\ws\tderive{\delta\isdq}$ (this changes only the linear combination $\lc$),
\barCl\label{eqn:obsnohf:2}
	\Tderive[2]{\isdq} &=& \frac{\np^2}{\Jl} v \transpose{\Equ{\isdq}}\J\delta\phisdq + M \J M^{-1}v\equ{\ws}\delta\w - \frac{\np}{\Jl}v\delta\Tl \ane \\
	& & +\>\lc(\delta\usdq, \delta\ws, \delta\isdq, \tderive{\delta\isdq}), \ans
\earCl
where $v:=\frac{\partial^2\Hmdq}{\partial{\phisdq}\partial{\phirdq}}\J\Equ{\phirdq}$. 

Gathering \eqref{eqn:obsnohf:is}, \eqref{eqn:obsnohf:1} and~\eqref{eqn:obsnohf:2}, we obtain the first six lines of the observability matrix
\beno
	\mathcal{O} := \begin{pmatrix}
		\frac{\partial^2\Hmdq}{\partial{\phisdq}^2} & \frac{\partial^2\Hmdq}{\partial{\phisdq}\partial{\phirdq}} & \Z{2}{1} & \Z{2}{1} \\
		-M & \Z{2}{2} & v & \Z{2}{1} \\
		\frac{\np^2}{\Jl} v \transpose{\Equ{\isdq}}\J & \Z{2}{2} &  M \J M^{-1}v\equ{\ws} & \frac{\np}{\Jl}v
	\end{pmatrix}.
\eeno
The first four lines of~$\mathcal{O}$ are independent, since the matrices $\frac{\partial^2\Hmdq}{\partial{\phisdq}\partial{\phirdq}}$ and~$M$ are invertible. Besides, as $M$~is not far from $\lambda\Id + \mu\J$, $M \J M^{-1}$~is not far from~$\J$; consequently the vectors $M \J M^{-1}v\equ{\ws}$ and~$\frac{\np}{\Jl}v$ are independent when~$\equ{\ws} \neq 0$. As a conclusion, the matrix~$\mathcal{O}$ is full rank when $\equ{\ws} \neq 0$, which means $\delta\phisdq$, $\delta\phirdq$, $\delta\w$ and $\delta\Tl$ can be obtained form $\isdq$, $\Tderive{\isdq}$ and $\Tderive[2]{\isdq}$.


On the contrary when~$\equ{\ws} = 0$, the last two lines of~$\mathcal{O}$ are collinear, as $v$~can be factored out. Besides, no extra information is gained by further time differentiating, since 
\barClno
	\Tderive[n]{\isdq} &=& \frac{\np^2}{\Jl}v\transpose{\Equ{\isdq}}\J \Tderive[n - 2]{\delta\phisdq} \\
		& & +\>\lc(\delta\usdq, \delta\ws, \delta\isdq, \tderive[n - 1]{\delta\isdq}) \\
	&=& \lc(\delta\usdq, \delta\ws, \delta\isdq, \tderive[n - 1]{\delta\isdq}),
\earClno
for~$n \geq 3$.

We have thus shown that the saturated IM is observable from the measured current~$\isdq$ when~$\equ{\ws} \neq 0$ and not observable when~$\equ{\ws} = 0$. This explains why sensorless control is more challenging at low speed. Notice the non-observability condition~$\equ{\ws} = 0$ corresponds to a line through the origin in the speed-torque plane.


		
\section{Analysis of HF signal injection}\label{sec:hf}
	\subsection{Virtual measurements and signal injection}\label{sec:effects}

A standard IM is geometrically non-salient, so it is paramount to understand magnetic saliency due to saturation to take advantage of signal injection. Following~\cite{Combes2016}, a precise understanding of saturation-induced saliency can be obtained, with:
\begin{itemize}
	\item a rigorous derivation of the ``HF model'' from the saturated model~\eqref{eqn:dynamic} and~\eqref{eqn:inv:currents}
	\item an analysis valid for any shape of injected signal (sinusoidal, square, \dots).
\end{itemize}

For the sake of simplicity, we restrict to pulsating signal injection as in~\cite{HaSIMS2000CRotIIAC}, namely
\be
	\usdq = \lf{\usdq} + \hf{\usdq}\finj(\W t),
\ee
where $\lf{\usdq}$~is the vector of control voltages, $\finj$~is a $1$\nbdash periodic function, $\W$~is a ``large'' constant, and $\hf{\usdq}$~is a vector determining the direction and amplitude of the injection.
The main result of~\cite{Combes2016} applied to the IM asserts that the state of~\eqref{eqn:dynamic} can be approximated by
\barCl\label{eqn:inj:state}
	\phisdq &=& \lf{\phisdq} + \frac{1}{\W}\hf{\usdq}\Finj(\W t) + \Landau\Bigl(\frac{1}{\W^2}\Bigr) \ase \label{eqn:inj:stator-flux} \\
	\phirdq &=& \lf{\phirdq} + \Landau\Bigl(\frac{1}{\W^2}\Bigr) \ase \label{eqn:inj:rotor-flux} \\
	\w &=& \lf{\w} + \Landau\Bigl(\frac{1}{\W^2}\Bigr) \ase \label{eqn:inj:speed}
\earCl
and the stator current measurements by
\be\label{eqn:inj:measurement}
	\isdq 
		  = \lf{\isdq} + \frac{1}{\W}\frac{\partial^2 \Hmdq}{\partial{\phisdq}^2}(\lf{\phisdq}, \lf{\phirdq})\hf{\usdq}\Finj(\W t) + \Landau\Bigl(\frac{1}{\W^2}\Bigr), 
\ee
where $\Finj$ is the primitive of $\finj$ with null mean, i.e.
\beno
	\Finj(\sigma) := \int_{0}^{\sigma} \finj(\tau) d\tau - \int_{0}^{1}\int_{0}^{\tau_1} \finj(\tau_2) d\tau_2 d\tau_1.
\eeno
In~\eqref{eqn:inj:state}--\eqref{eqn:inj:measurement}, $\lf{\phisdq}$, $\lf{\phirdq}$ and~$\lf{\w}$ are the state variables of the ``averaged system''
\barClno
	\Tderive{\lf{\phisdq}} &=& \lf{\usdq} - \Rs\lf{\isdq} - \J\ws\lf{\phisdq} \\
	\Tderive{\lf{\phirdq}} &=& -\Rr\lf{\irdq} - \J(\ws - \lf{\w})\lf{\phirdq} \\
	\frac{\Jl}{\np}\Tderive{\lf{\w}} &=& -\np \transpose{\lf{\phisdq}}\J\lf{\isdq} - \Tl,
\earClno
with the flux current-flux relations
\barClno
		\lf{\isdq} &=& \Pderive{\Hmdq}{\phisdq}(\lf{\phisdq}, \lf{\phirdq}) \\
		\lf{\irdq} &=& \Pderive{\Hmdq}{\phirdq}(\lf{\phisdq}, \lf{\phirdq}).
\earClno
By~\eqref{eqn:inj:measurement} the ``HF model'' is given by 
\be\label{eqn:inj:virtual}
	\hf{\isdq} := \frac{\partial^2 \Hmdq}{\partial{\phisdq}^2}(\lf{\phisdq}, \lf{\phirdq})\hf{\usdq},
\ee
$\frac{\hf{\isdq}}{\Omega}$~being the amplitude of the small current ripples created by the signal injection.
		The signals ~$\lf{\isdq}$ and~$\hf{\isdq}$ can then be recovered from the actual measurements~$\isdq$ by the filters
\barClno
	\lf{\isdq}(t) &\approx& \W \int_{t-\frac{1}{\W}}^{t} \isdq(\tau) d\tau \\
	\hf{\isdq}(t) &\approx& \W \frac{\int_{t-\frac{1}{\W}}^{t} \Bigl(\isdq\bigl(\tau - \frac{1}{2\W}\bigr) - \lf{\isdq}(\tau)\Bigr) \Finj(\W\tau) d\tau}{\int_{t-\frac{1}{\W}}^{t} \Finj^2(\W t) dt}. 
\earClno

The conclusion of this analysis is that signal injection yields the two extra ``virtual'' measurements~$\hf{\isdq}$ on top of the two ``physical'' measurements~$\lf{\isdq}$. We show in section~\ref{sec:obshf} that thanks to these extra measurements the observability can be recovered, which greatly simplifies the control problem.
	\subsection{The magnetic saliency matrix}\label{sec:xp}
		The essential part of the virtual measurement~\eqref{eqn:inj:virtual} is the 2x2 ``saliency matrix''
\beno
	\Sal(\lf{\phisdq}, \lf{\phirdq}) := \frac{\partial^2 \Hmdq}{\partial{\phisdq}^2}(\lf{\phisdq}, \lf{\phirdq}).
\eeno
It is moreover symmetric (it is the second derivative of a function),
hence can be parametrized as
\beno
\begin{pmatrix}
	a + b \cos{\sigma} & b \sin{\sigma} \\
	b \sin{\sigma} & a - b \cos{\sigma}
\end{pmatrix}
\eeno
where~$a$, $b \geq 0$ and~$\sigma \in ]-\pi, \pi]$ of course depend on~$\lf{\phisdq}$ and~$\lf{\phirdq}$.
With this parametrization the virtual measurements~\eqref{eqn:inj:virtual} read
\beno
\hf{\isdq} = a\hf{\usdq} + b\Bigl\|\hf{\usdq}\Bigr\|\begin{pmatrix}
	\cos(\sigma - \ti) \\
	\sin(\sigma - \ti)
\end{pmatrix},
\eeno
with $\ti$~the orientation of~$\hf{\usdq}$ in the $\dq$~frame.
The ``magnetic saliency'' is clearly visible in the sense that the relation between~$\hf{\usdq}$ and~$\hf{\isdq}$ depends on the absolute orientation~$\ti$.
Notice that in the unsaturated case, $a$~is a constant independent of~$\lf{\phisdq}$ and~$\lf{\phirdq}$ while $b$~is zero (this is easily seen from \eqref{eqn:model:linear}); there is no extra information in the virtual measurement~$\hf{\isdq}$, which is just $\hf{\usdq}$~times a constant. This is in accordance with the well-known fact that signal injection is useful for the control of the IM only because there is magnetic saturation.

\begin{figure}
	\setlength\figurewidth{6.5cm}
	\setlength\figureheight{2.5cm}
	\subfigure[$a$, the average saturation\label{fig:xp:a}]{\includegraphics{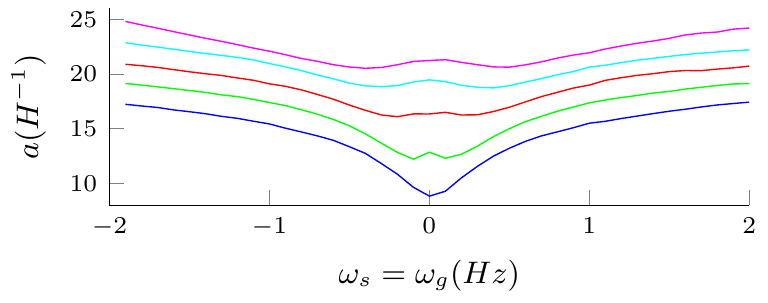}} \\
	\subfigure[$b$, the amplitude of the saliency\label{fig:xp:b}]{\includegraphics{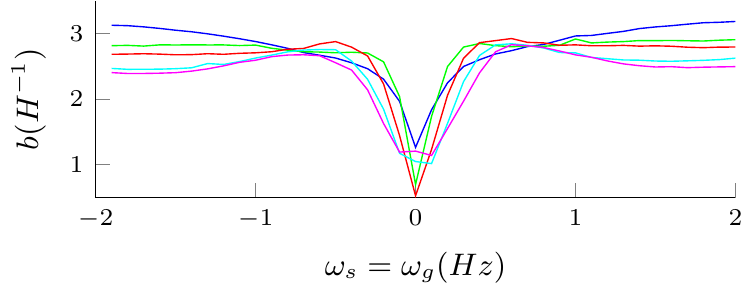}} \\
	\subfigure[$\sigma$, the orientation of the saliency\label{fig:xp:sigma}]{\includegraphics{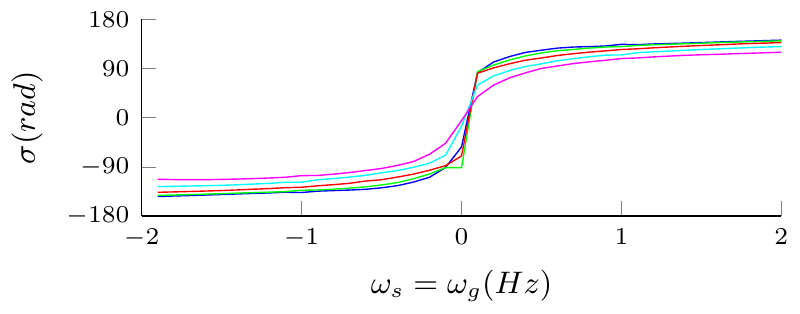}}
	\caption{Experimental characterization of magnetic saturation on a real $0.75kW$ IM under 
		$5\%$ (\textcolor{blue}{\bf---}), 
		$75\%$ (\textcolor{green}{\bf---}),
		$100\%$ (\textcolor{red}{\bf---}),
		$125\%$ (\textcolor{cyan}{\bf---}) and
		$150\%$ (\textcolor{magenta}{\bf---})
	of rated flux.}
	\label{fig:xp}
\end{figure}

We illustrate the shapes of~$a$, $b$ and~$\sigma$ with experimental data collected on a 0.75kW IM; the rated parameters listed in table~\ref{tbl:xp} corresponds to the usual unsaturated approximation around the rated point. Many tests with locked rotor were performed for various constant values of the stator speed~$\ws$ and currents~$\lf{\isdq}$. For each test at a given~$\ws$, we injected a square wave $\hf{\usdq}$ of amplitude 20V and frequency 500Hz with many different orientations~$\ti$; $\lf{\usdq}$~was chosen so that $\lf{\phirdq}$~is (approximately) at a prescribed level and aligned with the $d$~axis. We then reconstructed off-line the parameters~$a$, $b$ and~$\sigma$ corresponding to each experimental point from the recorded~$\hf{\isdq}$. The results are plotted in figure~\ref{fig:xp}.
%
%
%
It should be noted that the saliency orientation~$\sigma$ is not collinear with the direction of the rotor flux but almost orthogonal. To our knowledge, this has never been noticed even though \cite{YoonS2014ITOPE}~remarks that the saliency can slightly shift ($< 30^\circ$) from the rotor flux axis. 

\begin{figure}
	\setlength\figurewidth{6.5cm}
	\setlength\figureheight{2.5cm}
	\subfigure[$a$, the average saturation\label{fig:sim:a}]{\includegraphics{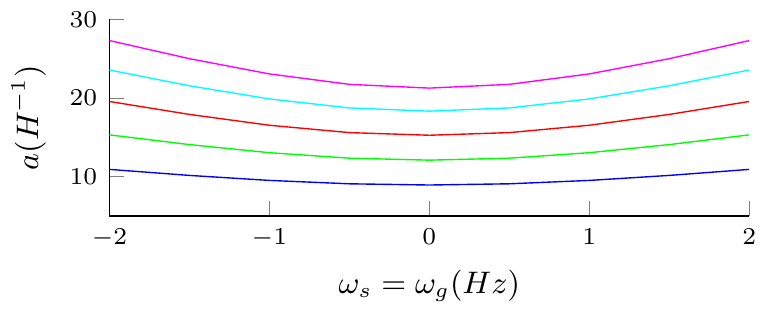}} \\
	\subfigure[$b$, the amplitude of the saliency\label{fig:sim:b}]{\includegraphics{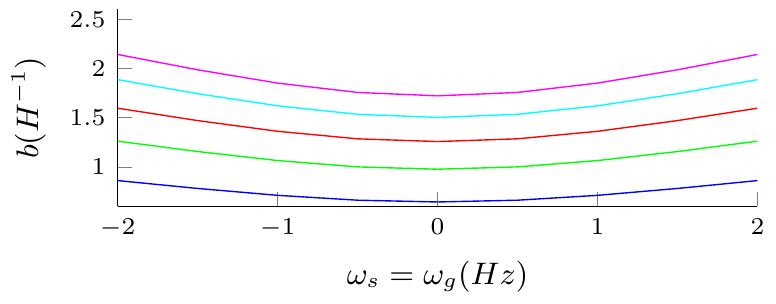}} \\
	\subfigure[$\sigma$, the orientation of the saliency\label{fig:sim:sigma}]{\includegraphics{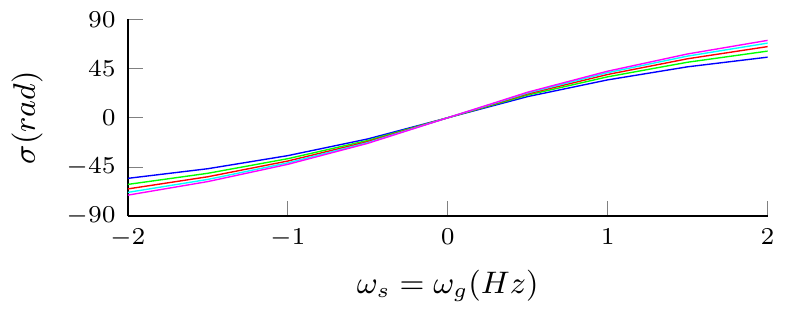}}
	\caption{Characterization of magnetic saturation on a simulated $0.75kW$ IM under 
		$5\%$ (\textcolor{blue}{\bf---}), 
		$75\%$ (\textcolor{green}{\bf---}),
		$100\%$ (\textcolor{red}{\bf---}),
		$125\%$ (\textcolor{cyan}{\bf---}) and
		$150\%$ (\textcolor{magenta}{\bf---})
		of rated flux.}
	\label{fig:sim}
\end{figure}

To fully take advantage of the approach, we would like to explicitly express the energy function~$\Hmdq$ from the experimental data. While this is rather easy for the PMSM (Permanent Magnet Synchronous Motor) \cite{JebaiMMR0IJoC}, it is much more difficult for the IM, since the rotor currents are not measured. We have not yet found a suitable parametrization of~$\Hmdq$, and therefore just illustrate the reasoning with the fictitious saturated energy function
\barCl\label{eqn:obs:model}
	\Hmdq(\phisdq, \phirdq) &=& \frac{1 + \epsm \norm{\phisdq + \phirdq}^2}{4(2\Lm + \Ll)}\norm{\phisdq + \phirdq}^2 \ane \\
		& & +\>\frac{1 + \epsl\norm{\phisdq + \phirdq}^2}{4\Ll}\norm{\phisdq - \phirdq}^2\!\!\!, \ans
\earCl
with the numerical parameters listed in table~\ref{tbl:sim}. This energy function is physically sensible, since it produces coefficients~$a$, $b$ and~$\sigma$ with reasonable magnitudes and shapes; it can be seen as a perturbation of the linear model~\eqref{eqn:model:linear}, and has of course the form~\eqref{eqn:model:sinus}. The data in figure~\ref{fig:sim} were obtained by performing in simulation the same procedure as in the experiments; the~$a$, $b$ and~$\sigma$ obtained by this ``simulated experiment'' are identical to the~$a$, $b$ and~$\sigma$ directly computed from the partial derivatives of~\eqref{eqn:obs:model}, which illustrates the relevance of the analysis in \ref{sec:effects}.

\begin{table}
	\centering
	\begin{tabular}{lr}
		\hline\noalign{\vskip 1pt}
		Rated power & $0.75kW$ \\
		Rated mechanical speed & $1500rpm$ \\
		Rated torque & $5N.m$ \\
		Rated voltage & $400V$ peak \\
		Rated current & $2A$ RMS \\
		Maximum current & $5A$ peak \\
		\hline\noalign{\vskip 1pt}
		Number of pole pairs & $2$ \\
		Inertia momentum & $5\cdot 10^{-3}kg.m^2$ \\
		Stator resistance & $13\Omega$ \\ 
		Rotor resistance & $10\Omega$ \\
		Mutual inductance & $0.42H$ \\
		Stator leakage inductance & $0.05H$ \\
		Rotor leakage inductance & $0.05H$ \\
		\hline
	\end{tabular}
	\caption{Rated motor characteristics and linear model parameters}
	\label{tbl:xp}
\end{table}

\begin{table}
	\centering
	\begin{tabular}{lr}
		\hline\noalign{\vskip 1pt}
		Mutual inductance $\Lm$ & $0.42H$ \\
		Leakage inductance $\Ll$ & $0.12H$ \\
		Mutual saturation factor $\epsm$ & $0.1Wb^{-2}$ \\
		Leakage saturation factor $\epsl$ & $1Wb^{-2}$ \\
		\hline
	\end{tabular}
	\caption{Magnetic parameters of the saturated IM model}
	\label{tbl:sim}
\end{table}

\section{Observability with signal injection}\label{sec:obshf}
		We now show that thanks to the virtual measurements~\eqref{eqn:inj:measurement}, the model becomes observable on the line~$\equ{\ws} = 0$. The linearized system~\eqref{eqn:model:linearized} can be written as
\be
	\Tderive{\delta x} = A \delta x + \begin{pmatrix}
		\Id \\
		\Z{3}{2}
	\end{pmatrix} \delta \usdq
\ee
where $\delta x := \transpose{(\transpose{\delta\phisdq}\!\!, \transpose{\delta\phirdq}\!\!,  \delta\w)}$ and
\beno
	A := \begin{pmatrix}
		-\Rs\frac{\partial \Hmdq}{\partial{\phisdq}^2} - \J \equ{\ws} & -\Rs\frac{\partial \Hmdq}{\partial{\phisdq}\partial{\phirdq}} & \Z{2}{1} \\
		-\Rr\frac{\partial \Hmdq}{\partial{\phirdq}\partial{\phisdq}} & -\Rr\frac{\partial \Hmdq}{\partial{\phirdq}^2} - \J\equ{\wg} & \J\phirdq \\
		\Z{1}{2} & \Z{1}{2} & 0
	\end{pmatrix}\!.
\eeno
The measured outputs of this linearized system are the linearized versions $\delta\lf{\isdq} = C\delta x$ and $\delta\hf{\isdq} = C_v\delta x$ of the actual measurement~$\lf{\isdq}$ and virtual measurement~$\hf{\isdq}$, where
\barClno
	C &:=& \begin{pmatrix}
		\frac{\partial^2\Hmdq}{\partial{\phisdq}^2} & \frac{\partial^2\Hmdq}{\partial{\phisdq}\partial{\phirdq}} & \Z{2}{1}
	\end{pmatrix} \\
	C_v &:=& \begin{pmatrix}
		\Pderive{}{\phisdq}\Bigl(\frac{\partial^2\Hmdq}{\partial{\phisdq}^2}\hf{\usdq}\Bigr) & 
		\Pderive{}{\phirdq}\Bigl(\frac{\partial^2\Hmdq}{\partial{\phisdq}^2}\hf{\usdq}\Bigr) & \Z{2}{1}
	\end{pmatrix}.
\earClno

\begin{figure}
	\setlength\figurewidth{6.5cm}
	\setlength\figureheight{3cm}
	\includegraphics{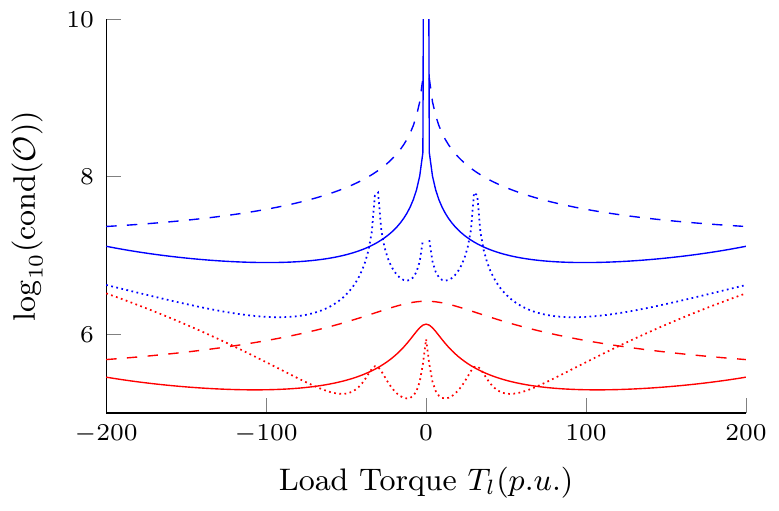}
	\caption{Condition numbers in logarithmic scale for observability matrices~$\mathcal{O}_s'$ (\textcolor{blue}{\bf---})
	and~$\mathcal{O}_s$ (\textcolor{red}{\bf---}) at~$\equ{\ws} = 0$ and~$50\%$ (dotted), $100\%$ (solid) and~$150\%$ (dashed) of nominal flux.}
	\label{fig:cond}
\end{figure}

Expressing $\delta\phisdq$, $\delta\phirdq$ and $\delta\w$ in function of $\delta\lf{\isdq}$, $\Tderive{\delta\lf{\isdq}}$, $\delta\hf{\isdq}$ and $\Tderive{\delta\hf{\isdq}}$ is possible if and only if the observability matrix 
\beno
	\mathcal{O}_s := \begin{pmatrix}
		C \\
		C_v \\
		C A \\
		C_v A
	\end{pmatrix}
\eeno
has full rank 5. Differentiating $\delta\lf{\isdq}$ once more then yields $\delta\Tl$. Analytically checking that $\mathcal{O}_s$ has full rank is very tedious; therefore we will make do with numerical computations based on the energy function \eqref{eqn:obs:model}, which is representative of the general situation. Moreover, rather than computing the rank of $\mathcal{O}_s$, we compute its condition number (i.e. the ratio between its largest and smallest singular values). Indeed, a large condition number means the matrix is ``nearly not full rank'', and make numerical inversion very inaccurate in practice. Figure~\ref{fig:cond} shows plots of this condition number in function of the load torque for different values of the rotor flux on the unobservability line $\equ{\ws} = 0$; we can see that first-order observability is restored thanks to signal injection (recall that without signal injection, the condition number is infinite since the system is not observable). It is therefore possible to control the IM at zero stator velocity with a basic control law based for instance on an extended Kalman filter.

Notice that the control laws of~\cite{Brandstetter2011,JanseL1995IToIA,YoonS2014ITOPE,Zatocil2008} do not use an observer to reconstruct the stator and rotor fluxes, but algebraically compute them from~$\lf{\isdq}$ and~$\hf{\isdq}$. The rotor velocity is then obtained from a PI loop considering the fluxes are known. Assuming the rotor flux is aligned with the d axis, this corresponds to inverting the matrix
\beno
	\mathcal{O}_s' := \begin{pmatrix}
		C \\
		C_v \\
		(C A)_2
	\end{pmatrix}\!,
\eeno
where $(C A)_2$~denotes the second line of~$C A$. But the condition number of~$\mathcal{O}_s'$ may be very large in some operating conditions (see figure~\ref{fig:cond}). This means the control laws~\cite{Brandstetter2011,JanseL1995IToIA,YoonS2014ITOPE,Zatocil2008} may perform poorly in those conditions.
\section{Conclusion}
		We have shown that thanks to magnetic saturation and signal injection, first-order observability is restored at any operating point, in particular on the unobservability line $\equ{\ws} = 0$. This means that low-speed sensorless control of the IM is doable even with model uncertainties, since it is possible to reconstruct the state with an essentially linear observer (fixed-gain observer, gain-scheduled observer, or extented Kalman filter). 


\bibliographystyle{phmIEEEtran}
\bibliography{biblio}

%
%
%
%

\end{document}